%&amstex          
\input amstex\documentstyle{amsppt}  
\pagewidth{12.5cm}\pageheight{19cm}\magnification\magstep1
\topmatter
\title{On bases of certain Grothendieck groups, II}\endtitle
\author G. Lusztig\endauthor
\address{Department of Mathematics, M.I.T., Cambridge, MA 02139}\endaddress
\thanks{Supported by NSF grant DMS-2153741}\endthanks
\endtopmatter   
\document

\define\si{\sim}

\define\sqc{\sqcup}

\define\lb{\linebreak}

\define\op{\oplus}
   
\define\part{\partial}
\define\emp{\emptyset}

\define\n{\notin}

\define\m{\mapsto}
\define\do{\dots}

\define\sub{\subset}

\define\ti{\tilde}
\define\nl{\newline}
\redefine\i{^{-1}}

\define\un{\underline}

\define\supp{\text{\rm supp}}

\define\a{\alpha}

\define\g{\gamma}

\define\e{\epsilon}

\define\io{\iota}

\define\p{\pi}

\define\ps{\psi}

\redefine\t{\tau}

\define\z{\zeta}

\define\Ph{\Phi}

\define\QQ{\bold Q}

\define\ZZ{\bold Z}

\define\cp{\Cal P}

\define\cs{\Cal S}

\define\cx{\Cal X}

\define\tE{\ti E}

\define\tY{\ti Y}

\define\sha{\sharp}

\define\bE{\bar E}

\head Introduction\endhead
\subhead 0.1\endsubhead
Let $\cs$ be the set of (isomorphism classes of) unipotent
representations of

(a) a symplectic or odd special orthogonal group over a finite field, or

(b) an even split special orthogonal group over a finite field.
\nl
Let $W$ be the corresponding Weyl group and let $Ce(W)$ be the set of
two-sided cells of $W$. It is known that we have a natural partition
$\cs=\sqc_{c\in Ce(W)}\cs_c$. Thus the Grothendieck group $\ZZ[\cs]$ is
a direct sum $\op_{c\in Ce(W)}\ZZ[\cs_c]$ where for a finite set $V$,
$\ZZ[V]$ denotes the free abelian group with basis $V$. We now fix
$c\in Ce(W)$. In \cite{L20} a new basis of $\ZZ[\cs_c]$ was defined
in case (a)
with strong positivity properties with respect to Fourier transform.
In case (b) there were two versions of such a basis in \cite{L20};
later the second version was adopted. (It is the only version
compatible with the general scheme of \cite{L23}.)
However, this second version
was given in \cite{L20} without proof. It was further studied
in \cite{L22} but again with several statements given without proof.
The purpose of this paper is to provide proofs of the
statements whose proofs were missing in \cite{L20}, \cite{L22}.
One of the ideas in \cite{L22} was to enlarge the setup of (b) to
include also the unipotent representations of a full orthogonal group
both in the split and in the nonsplit case. In this enlarged setup
it is possible to make a comparison with case (a) and to deduce all
desired statements from the known ones in case (a).

\subhead 0.2. Notation\endsubhead
Let $\ti0=2\ZZ,\ti1=2\ZZ+1$.
For $n\in\ZZ$ we define $\un n\in\ZZ$ by $\un n=0$ if $n\in\ti0$,
$\un n=1$ if $n\in\ti1$. Let $F=\ZZ/2\ZZ$.
The number of elements in a finite set $X$ is denoted by $\sha X$.
For $i\le j$ in $\ZZ$ we set $[i,j]=\{h\in\ZZ;i\le h\le j\}$.

\head 1. The bijection $\z:\cx_{N-2}@>\si>>\cx_{N-1}$\endhead
\subhead 1.1 \endsubhead
We fix $N\in\{3,5,7,\do\}$.
 Let $\tE_N$ be the set of all subsets of $S_N:=[1,N]$,
 an $F$-vector space of dimension $N$ in which the sum of $X,X'$ is
$(X\cup X')-(X\cap X')$.

 Let $E_N=\{X\in\tE_N;\sha X\in\ti0\}$, an $F$-subspace of $\tE_N$ of
dimension $N-1$.
Let ${}^2E_N=\{X\in E_N;\sha X=2\}$.
A subset $\{i,j\}\in{}^2E_N$ will be often written as $ij$ if
$i<j$ and $j-i\in\ti1$ or if $i>j$ and $i-j\in\ti0$.

Let $\cp_N$ be the set of all unordered sequences $X_1,X_2,\do,X_k$ in
${}^2E_N$, $(k\ge0)$ such that for any $a\ne b$ we have
$X_a\cap X_b=\emp$.

\subhead 1.2\endsubhead
For $B\in\cp_N$ let $\supp(B)=\cup_{\{i,j\}\in B}\{i,j\}\sub S_N$.
For $B\in\cp_N$ we have $B=B^0\sqc B^1$ where
$B^0=\{\{i,j\}\in B;i-j\in\ti0\}$, $B^1=\{\{i,j\}\in B;i-j\in\ti1\}$.
Let $B\in\cp_N$ and let $X\in\tE_N$. We
say that $X$ is $0$-covered (resp. $1$-covered) by $B^1$ if there exist
$a_1b_1\in B^1,a_2b_2\in B^1,\do,a_sb_s\in B^1$,
such that

$X=[a_1,b_1]\sqc[a_2,b_2]\sqc[a_s,b_s]$

(resp. $X=[a_1,b_1]\sqc[a_2,b_2]\sqc[a_s,b_s]\sqc\{u\}$ for some $u$).
\nl
Let $B\in\cp_N$. We say that $B\in{}^*\cp_N$ if the following
conditions are satisfied:

For every $ij\in B^1$, $[i+1,j-1]$ is $0$-covered by $B^1$. There
exists a sequence $i_*(B)=(i_1,i_2,i_3,\do,i_{2s})$ in $S_N$ such that
$B^0=\{i_{2s}i_1,i_{2s-1}i_2,\do,i_{s+1}i_s\}$ (it is automatically
unique and we have $s=\sha B^0$). If $s\ge1$, each of the intervals

$[i_1+1,i_2-1],[i_2+1,i_3-1],\do,[i_{s-1}+1,i_s-1]$,

$[i_{s+1}+1,i_{s+2}-1],\do,[i_{2s-1}+1,i_{2s}-1]$
\nl
is $0$-covered by $B^1$. 

(It follows that any two consecutive terms of $i_1,i_2,\do,i_s$ have
different parities and any two consecutive terms of
$i_{s+1},i_{s+2},\do,i_{2s}$ have different parities. Also, $i_1,i_{2s}$' have the same parity.)

\subhead 1.3\endsubhead
Let $B\in{}^*\cp_N$ and let
$i_*(B)=(i_1,\do,i_{2s})$.
We consider the following conditions on $B$.

(I) We have either $s=0$ or $s\ge1$ and each of $[1,i_1-1],[i_{2s}+1,N]$
is $0$-covered by $B^1$ (in particular, $i_1,i_{2s}$ are odd).

(II) We have $N\n\supp(B)$ and either $s=0$ or else $s$ is odd and one of
(i),(ii) below holds:

(i) $[1,i_1-1]$ is $1$-covered by $B^1$ and $[i_{2s}+1,N-1]$ is
$0$-covered by $B^1$. 

(ii) $[1,i_1-1]$ is $0$-covered by $B^1$ and $[i_{2s}+1,N-1]$ is
$1$-covered by $B^1$. 

(III)  Same as in (I). In addition, if $s$ is even, then $\{i,N\}\in B$
for some even $i$; if $s$ is odd, then $\{i,N\}\in B$ for some odd $i$.

Let

$\cx_{N-1}=\{B\in{}^*\cp_N; B \text{ satisfies } (I)\}$,

$\cx_{N-2}^+=\{B\in{}^*\cp_N; B \text{ satisfies } (II)\}$,

$\cx_{N-2}^-=\{B\in{}^*\cp_N; B \text{ satisfies } (III)\}$, 

$\cx_{N-2}=\cx_{N-2}^+\sqc\cx_{N-2}^-$.

Let $\ti\cx_{N-1},\ti\cx_{N-2}$ be the subsets of ${}^*\cp_N$ defined
in \cite{L22,1.6}. It is immediate that $\cx_{N-1}=\ti\cx_{N-1}$,
$\cx_{N-2}=\ti\cx_{N-2}$. Using \cite{L22,1.6,1.7}, we see also that
$\cx_{N-1}$ (resp. $\cx_{N-2}$) coincides with what in \cite{L22} was
denoted by $\cx_{N-1}$ (resp. $\cx_{N-2}$).

Using \cite{L22,1.8} we see that for $B\in\cx_{N-2}^+$ in case (II)(i)
there is a unique $u=u_B\in[1,i_1-1]$ such that $[1,u-1]$ and
$[u+1,i_1-1]$ are $0$-covered by $B^1$ (note that $i_1\in\ti0,u_B\in\ti1$),
while for $B\in\cx_{N-2}^+$ in case (II)(ii) there is a unique
$u=u_B\in[i_{2s}+1,N-1]$ such that $[i_{2s}+1,u-1]$ and $[u+1,N-1]$ are
$0$-covered by $B^1$ (note that $i_{2s}\in\ti1,u_B\in\ti0$).

\subhead 1.4\endsubhead
Following \cite{L22,1.9} for any $t\in\ti0$ we set

$\cx_{N-1}^t=\{B\in\cx_{N-1};\sha B^0=t\}$ if $t\ge0$,

$\cx_{N-1}^t=\{B\in\cx_{N-1};\sha B^0=-t-1\}$ if $t<0$,

$\cx_{N-2}^{t,+}=\{B\in\cx_{N-2}^+;\sha B^0=0\}$ if $t=0$,

$\cx_{N-2}^{t,+}=\{B\in\cx_{N-2}^+;\sha B^0=t-1,u_B\in\ti1\}$ if
$t\in\{2,4,6,\do\}$,

$\cx_{N-2}^{t,+}=\{B\in\cx_{N-2}^+;\sha B^0=-t-1,u_B\in\ti0\}$ if
$t\in\{-2,-4,-6,\do\}$,

$\cx_{N-2}^{t,-}=\{B\in\cx_{N-2}^-;\sha B^0=0\}$ if $t=0$,

$\cx_{N-2}^{t,-}=\{B\in\cx_{N-2}^-;\sha B^0=t,\{i,N\}\in B
\text{ for some even }i\}$ if $t\in\{2,4,6,\do\}$,

$\cx_{N-2}^{t,-}=\{B\in\cx_{N-2}^-;\sha B^0=-t-1,\{i,N\}\in B
\text{ for some odd }i\}$

if $t\in\{-2,-4,-6,\do\}.$

We have

$\cx_{N-1}=\sqc_{t\in\ti0}\cx_{N-1}^t$,

$\cx_{N-2}^+=\sqc_{t\in\ti0}\cx_{N-2}^{t,+}$,

$\cx_{N-2}^-=\sqc_{t\in\ti0}\cx_{N-2}^{t,-}$.

\subhead 1.5\endsubhead
For any even $s\ge0$ we define a bijection

$$\{B\in\cx_{N-2}^-;\sha B^0=s\}@>\si>>
\{B\in\cx_{N-1};\sha B^0=s,\{i,N\}\in B \text{ for some even }i\}\tag a$$

by $B\m B$. For any odd $s>0$ we define a bijection
$$\{B\in\cx_{N-2}^-;\sha B^0=s\}@>\si>>
\{B\in\cx_{N-1};\sha B^0=s,\{i,N\}\in B \text{ for some odd } i\}\tag b$$

by $B\m B$. We define a bijection
$$\{B\in\cx_{N-2}^+;\sha B^0=0\}@>\si>>
\{B\in\cx_{N-1};\sha B^0=0,N\n\supp(B)\}\tag c$$

by $B\m B$. For any odd $s>0$ we define a bijection

$$\align&\{B\in\cx_{N-2}^+;\sha B^0=s,u_B\in\ti1\}@>\si>>\\&
\{B\in\cx_{N-1};
\sha B^0=s+1,\{i,N\}\in B \text{ for some odd }i\}\tag d
\endalign$$

by $B\m B\sqc\{u_B,N\}$. For any odd $s>0$ we define a bijection

$$\align&\{B\in\cx_{N-2}^+;\sha B^0=s,u_B\in\ti0\}@>\si>>\\&
\{B\in\cx_{N-1};
\sha B^0=s,\{i,N\}\in B \text{ for some even }i\}\tag e\endalign$$

by $B\m B\sqc\{u_B,N\}$.

Putting together these bijections we get a bijection

(f) $\z:\cx_{N-2}@>\si>>\cx_{N-1}$
\nl
which restricts to bijections

(g) $\cx_{N-2}^-@>\si>>
\{B\in\cx_{N-1};\{i,N\}\in B \text{ for some $i$ such that }
i=\sha B^0\mod2\}$,

(h) $\cx_{N-2}^+@>\si>>
\{B\in\cx_{N-1};\{i,N\}\n B \text{ for any $i$ such that }
i=\sha B^0\mod2\}$.

Let $\cx_{N-1,0}$ be the set of all $B\in\cx_{N-1}$ such that
$\sha B^0>0$ and such that $\{i,N\}\in B$ for a unique $i=i^B$
with $i=\sha B^0+1\mod2$.

The inverse bijection $\z\i:\cx_{N-1}@>\si>>\cx_{N-2}$ takes any
$B\in\cx_{N-1,0}$ to $B-\{i^B,N\}$ and takes any other $B$ to $B$.

\head 2. The
bijections $\e:\cx_{N-1}@>\si>>E_N,{}'\e:\cx_{N-2}@>\si>>E_N$
\endhead
\subhead 2.1\endsubhead
For $X\in{}^2E_N$ we define $[[X]]\in E_N$ as follows:

$[[X]]=\{i,i+1,\do,j\}=[i,j]$

if $X=ij$ with $i<j,j-i\in\ti1$;

$[[X]]=\{i,i+1,\do,N-1,N,1,2,\do,j\}=[i,N]\cup[1,j]$

if $X=ij$ with $i>j,i-j\in\ti0$.

\subhead 2.2\endsubhead
Let $B\in\cx_{N-1}$. Following \cite{L22} we set:

(a) $\e(B)=\sum_{ij\in B}[[ij]]\in E_N$.
\nl
(This definition is a reformulation of a definition which appears
in \cite{L20}.)

\subhead 2.3\endsubhead
For $B=\{X_1,X_2,\do,X_k\}\in\cp_N$ let $<B>$ be the $F$-subspace of
$E_N$ spanned by $X_1,X_2,\do,X_k$. The following result is proved in
\cite{L20}.

\proclaim{Theorem 2.4} (a) The map $\cx_{N-1}@>>>E_N$ given by $B\m\e(B)$
is a bijection.

(b) If $B\in\cx_{N-1}$, then $\e(B)\in<B>$.

(c) If $B\in\cx_{N-1},B'\in\cx_{N-1}$ are such that $\e(B')\in<B>$,
then we have $\sha B'{}^0\le\sha B^0$.
\endproclaim

\subhead 2.5\endsubhead
For $B,B'$ in $\cx_{N-1}$ we write $B'\le B$ if there exists a sequence

$B'=B_0,B_1,\do,B_h=B$
\nl
in $\cx_{N-1}$ such that 
$$\e(B_0)\sub<B_1>,\e(B_1)\sub<B_2>,\do,\e(B_{h-1})\sub<B_h>.$$

The following result is proved in \cite{L20}.

(a) {\it $\le$ is a partial order on $\cx_{N-1}$.}
\nl
Using 2.4(c) we see that:

(b) {\it If $B,B'$ in $\cx_{N-1}$ satisfy $B'\le B$, then
$\sha B'{}^0\le\sha B^0$.}

\subhead 2.6\endsubhead
For $X\in E_N$ we write $<X>=<B>$ where $B=\e\i(X)$; this is an
$F$-subspace of $E_N$.
For $X,X'$ in $E_N$ we write $X'\le X$ if $\e\i(X')\le\e\i(X)$. By
2.5(a) this is a partial order on $E_N$. Let $\ZZ[E_N]$ be the abelian
group consisting of functions $E_N@>>>\ZZ$. This has a $\ZZ$-basis
$\{\ps_X;X\in E_N\}$ where for $X\in E_N$, $\ps_X:E_N@>>>\ZZ$ is given by
$X\m1,X'\m0$ if $X'\in E_N-\{X\}$. For any subset $Z$ of $E_N$ we set
$\ps_Z=\sum_{X'\in Z}\ps_{X'}\in\ZZ[E_N]$. Thus for $X\in E_N$ we have
$\ps_X=\ps_{\{X\}}$.
Note that $\{\ps_{<X>};X\in E_N\}$ is related to
$\{\ps_X;X\in E_N\}$ 
by an upper triangular matrix with entries in $\ZZ$ and with
$1$ on diagonal.
In particular, $\{\ps_{<X>};X\in E_N\}$ is a $\ZZ$-basis of
$\ZZ[E_N]$.

\subhead 2.7\endsubhead
For $X\in E_N$ we set $\g(X)=\sha(X\cap\ti0)-\sha(X\cap\ti1)\in\ti0$.
For $t\in\ti0$ we set $E_N^t=\{X\in E_N;\g(X)=t\}$. We have
$E_N=\sqc_{t\in\ti0}E_N^t$. In \cite{L22,1.15} it is shown that for any 
$t\in\ti0$, the map $\e:\cx_{N-1}@>>>E_N$ takes $\cx_{N-1}^t$ into
$E_N^t$. Combining this with 2.4(a), we see that
$\e$ induces for any $t\in\ti0$ a bijection $\cx_{N-1}^t@>\si>>E_N^t$.

\subhead 2.8\endsubhead
For $B\in\cx_{N-2}$ we define ${}'\e(x)\in E_N$ as follows.

If $B\in\cx_{N-2}^+$ satisfies $\sha B^0=0$, we set
${}'\e(B)=\sum_{ij\in B}[[ij]]$.

If $B\in\cx_{N-2}^+$ satisfies $\sha B^0\in\ti1,u_B\in\ti0$, we set
${}'\e(B)=\sum_{ij\in B}[[ij]]+[u_B,N]$.

If $B\in\cx_{N-2}^+$ satisfies $\sha B^0\in\ti1,u_B\in\ti1$, we set

${}'\e(B)=\sum_{ij\in B}[[ij]]+[N,N]+[1,u_B]$.

If $B\in\cx_{N-2}^-$ we set ${}'\e(B)=\sum_{ij\in B}[[ij]]$.

From the definitions, for $B\in\cx_{N-2}$ we have

(a) ${}'\e(B)=\e(\z(B))$.

Part (a) of the following result has been stated without proof in
\cite{L22}.

\proclaim{Theorem 2.9} (a) The map ${}'\e:\cx_{N-2}@>>>E_N$ is a bijection.

(b) If $B\in\cx_{N-2}$ then ${}'\e(B)\in<B>$.
\endproclaim
Indeed, the map in (a) is the composition of the bijection $\z$ with the
bijection $\e$ in 2.4(a). Part (b) is proved in \cite{L22,1.7}.

\subhead 2.10 \endsubhead
For $B,B'$ in $\cx_{N-2}$ we wite $B'\preceq B$ if there exists a
sequence

$B'=B_0,B_1,\do,B_h=B$
\nl
in $\cx_{N-2}$ such that 
$${}'\e(B_0)\sub<B_1>,{}'\e(B_1)\sub<B_2>,\do,{}'\e(B_{h-1})
\sub<B_h>.\tag a$$

The following result was stated without proof in \cite{L22}.

(b) {\it $\preceq$ is a partial order on $\cx_{N-2}$.}
\nl
For any $\ti B\in\cx_{N-2}$ we have (using the definition of $\z$)
$\ti B\sub\z(\ti B)$. Hence if (a) holds then for $k=0,1,\do,h-1$ we
have $\e(\z(B_k))\sub<B_k>\sub<\z(B_{k+1})>$, so that
$\e(\z(B_k))\sub<\z(B_{k+1})>$. We see that if $B'\preceq B$,  then
$\z(B')\le\z(B)$. It is enough to show that if $B'\preceq B$ and
$B\preceq B'$ then $B=B'$.
But our assumption implies $\z(B')\le\z(B)$ and $\z(B)\le\z(B')$
hence by 2.5(a) we have $\z(B)=\z(B')$. Since $\z$ is a bijection
we deduce that $B=B'$, proving (b).

The proof above shows that:

(c) {\it For $B',B$ in $\cx_{N-2}$ such that $B'\preceq B$ we have
$\z(B')\le\z(B)$. }

\subhead 2.11\endsubhead
For $X\in E_N$ we write ${}'<X>=<B>$ where $B={}'\e\i(X)$; this is an
$F$-subspace of $E_N$. For $X,X'$ in $E_N$ we write $X'\preceq X$ if
${}'\e\i(X')\preceq{}'\e\i(X)$. By 2.10(b) this is a partial order on
$E_N$.
Note that $\{\ps_{{}'<X>};X\in E_N\}$ is related to
$\{\ps_X;X\in E_N\}$ 
by an upper triangular matrix with entries in $\ZZ$ and with
$1$ on diagonal.
In particular, $\{\ps_{{}'<X>};X\in E_N\}$ is a $\ZZ$-basis of
$\ZZ[E_N]$.

\subhead 2.12\endsubhead
We have a partition $E_N=E_N^+\sqc E_N^-$ where
$$E_N^+=\{X\in E_N;N\n X\},E_N^-=\{X\in E_N;N\in X\}.$$
For $t\in\ti0$ we set

$E_N^{t,+}=E_N^t\cap E_N^+,E_N^{t,-}=E_N^t\cap E_N^-$.
\nl
We have

$E_N^+=\sqc_{t\in\ti0}E_N^{t,+}$, $E_N^-=\sqc_{t\in\ti0}E_N^{t,-}$.
\nl
In \cite{L22,1.15} it is shown that for any 
$t\in\ti0$ the map ${}'\e:\cx_{N-2}@>>>E_N$ takes $\cx_{N-2}^{t,+}$ into
$E_N^{t,+}$ and $\cx_{N-2}^{t,-}$ into $E_N^{t,-}$; in particular it
takes $\cx_{N-2}^+$ into $E_N^+$ and $\cx_{N-2}^-$ into $E_N^-$.
Combining this with 2.9(a) we obtain the following result which was
stated without proof in \cite{L22,1.15}.

\proclaim{Theorem 2.13} ${}'\e:\cx_{N-2}@>>>E_N$ induces bijections
$\cx_{N-2}^+@>\si>>E_N^+$ and $\cx_{N-2}^-@>\si>>E_N^-$; moreover,
for any $t\in\ti0$ it induces a bijection
$\cx_{N-2}^{t,+}@>\si>>E_N^{t,+}$ and a bijection
$\cx_{N-2}^{t,-}@>\si>>E_N^{t,-}$. 
\endproclaim

\subhead 2.14\endsubhead
Let $B',B$ in $\cx_{N-2}$ be such that $B'\preceq B$. We show:

(a) {\it If $B\in\cx_{N-2}^+$ then $B'\in\cx_{N-2}^+$.}
\nl
We can assume that ${}'\e(B')\in<B>$. (The general case follows by
applying this special case several times.) Assume that
$B'\in\cx_{N-2}^-$. By 2.13 we then have  $N\in{}'\e(B')$. We have
$N\n\supp(B)$. Hence for any $X\in<B>$ we have $N\n X$. Since
${}'\e(B')\in<B>$, we have $N\n{}'\e(B')$. This is a contradiction; (a)
is proved.

\subhead 2.15\endsubhead
Let $\ZZ[E_N^+]$ (resp. $\ZZ[E_N^-]$) be the subgroup of $\ZZ[E_N]$ with
$\ZZ$-basis $\{\ps_X;X\in E_N^+\}$ (resp. $\{\ps_X;X\in E_N^-\}$).
By 2.14(a), for $X\in E_N^+$ we have $\ps_{{}'<X>}\in\ZZ[E_N^+]$.
Note that $\{\ps_{{}'<X>};X\in E_N^+\}$ is related to
$\{\ps_X;X\in E_N^+\}$ by an upper triangular matrix
with entries in $\ZZ$ and with $1$ on diagonal.
Hence $\{\ps_{{}'<X>};X\in E_N^+\}$ is a $\ZZ$-basis of $\ZZ[E_N^+]$.
Also, $\{\ps_{{}'<X>\cap E_N^-};X\in E_N^-\}$ is related to
$\{\ps_X;X\in E_N^-\}$
by an upper triangular matrix with entries in $\ZZ$ and with $1$ on diagonal.
Hence $\{\ps_{{}'<X>\cap E_N^-};X\in E_N^-\}$ is a $\ZZ$-basis of
$\ZZ[E_N^-]$.

\head 3. The partial order $\bar\preceq$\endhead
\subhead 3.1\endsubhead
Let

${}'\cx_{N-2}^{0,+}=\{B\in\cx_{N-2}^{0,+};N-1\n\supp(B)\}$,

${}''\cx_{N-2}^{0,+}=\{B\in\cx_{N-2}^{0,+};N-1\in\supp(B)\}$,

${}'E_N^{0,+}=\{X\in E_N^{0,+}; N-1\n X\}$,

${}''E_N^{0,+}=\{X\in E_N^{0,+}; N-1\in X\}$.
\nl
We have

$\cx_{N-2}^{0,+}={}'\cx_{N-2}^{0,+}\sqc{}''\cx_{N-2}^{0,+}$,

$E_N^{0,+}={}'E_N^{0,+}\sqc{}''E_N^{0,+}$.
\nl
By \cite{L22,4.3}, the map ${}'\e:\cx_{N-2}^{0,+}@>>>E_N^{0,+}$
    restricts to maps
$${}'\cx_{N-2}^{0,+}@>>>{}'E_N^{0,+},
{}''\cx_{N-2}^{0,+}@>>>{}''E_N^{0,+}.$$

Using 2.9(a), it follows that

(a) {\it these two maps are bijections.}
\nl
(This was stated without proof in \cite{L22}.)

We set

${}'\cx_{N-2}^+
={}'\cx_{N-2}^{0,+}\sqc\sqc_{t\in\ti0,t<0}\cx_{N-2}^{t,+}$,

${}''\cx_{N-2}^+
={}''\cx_{N-2}^{0,+}\sqc\sqc_{t\in\ti0,t>0}\cx_{N-2}^{t,+}$,

${}'\cx_{N-2}^-=\sqc_{t\in\ti0,t\ge0}\cx_{N-2}^{t,-}$,

${}''\cx_{N-2}^-=\sqc_{t\in\ti0,t<0}\cx_{N-2}^{t,-}$,

${}'E_N^+={}'E_N^{0,+}\sqc\sqc_{t\in\ti0,t<0}E_N^{t,+}$,

${}''E_N^+={}''E_N^{0,+}\sqc\sqc_{t\in\ti0,t>0}E_N^{t,+}$,

${}'E_N^-=\sqc_{t\in\ti0,t\ge0}E_N^{t,-}$,

${}''E_N^-=\sqc_{t\in\ti0,t<0}E_N^{t,-}.$

Note that ${}'\e$ defines bijections

${}'\cx_{N-2}^+@>\si>>{}'E_N^+$,

${}''\cx_{N-2}^+@>\si>>{}''E_N^+$,

${}'\cx_{N-2}^-@>\si>>{}'E_N^-$,

${}''\cx_{N-2}^-@>\si>>{}''E_N^-$.

\subhead 3.2\endsubhead
Let $B'\in\cx_{N-2}^{t',+},B\in\cx_{N-2}^{t,+}$ be such that
$B'\preceq B$; here $t\in\ti0,t'\in\ti0$. The following result was
stated without proof in \cite{L22,4.4}.

(a) {\it If $t<0$ then we have $t=t'$ or $|t'|<|t|$.   }
\nl
We have $\z(B')\le\z(B)$, see 2.10(c);
hence $\sha \z(B')^0\le\sha\z(B)^0$, see 2.5(b). If $t'<0$, then
$\sha\z(B)^0=\sha B^0=-t-1$,
$\sha\z(B')^0=\sha B'{}^0=-t'-1$, so that $-t'-1\le-t-1$; thus $t'=t$ or
$-t'<-t$. If $t'>0$, then $\sha\z(B)^0=\sha B^0=-t-1$,
$\sha\z(B')^0=\sha B'{}^0+1=t'$, so that $t'\le-t-1<-t$ and
$|t'|<|t|$.
If $t'=0$ the result is obvious. This proves (a).

\subhead 3.3\endsubhead
Let $B'\in\cx_{N-2}^+,B\in{}'\cx_{N-2}^{0,+}$ be such that $B'\preceq B$.
The following result was stated without proof in \cite{L22,4.4}.

(a) {\it We have $B'\in{}'\cx_{N-2}^{0,+}$.}
\nl
We can assume that ${}'\e(B')\in<B>$. (The general case follows by
applying this special case several times.) 
We have $B'\in{}'\cx_{N-2}^{t',+}$ for some $t'\in0$.
We have $\z(B')\le\z(B)$, see 2.10(c); hence
$\sha\z(B')^0\le\sha\z(B)^0$, see 2.5(b).
If $t'>0$, then $\sha\z(B)^0=\sha B^0=0$, $\sha\z(B')^0=
\sha B'{}^0+1=t'$, so
that $t'\le0$, contradiction.
If $t'<0$, then $\sha \z(B)^0=\sha B^0=0$, $\sha\z(B')^0=
\sha B'{}^0=-t'-1$, so
that $-t'-1\le0$ and $t'\ge-1$. Thus $t'\le-1$ and $t'\ge-1$ so that
$t'=-1$, contradicting $t'\in\ti0$. We see that $t'=0$.
Assume that $B'\in{}''\cx_{N-2}^{0,+}$. Then by 3.1(a) we have
${}'\e(B')\in{}''E_N^{0,+}$ so that $N-1\in{}'\e(B')$. 
We have $N-1\n\supp(B)$. Hence for any $X\in<B>$ we have $N-1\n X$. Since
${}'\e(B')\in<B>$, we have $N-1\n{}'\e(B')$. This is a contradiction.
We see that $B'\in{}'\cx_{N-2}^{0,+}$; (a) is proved.

\subhead 3.4\endsubhead
Let $B'\in\cx_{N-2}^{t',-},B\in\cx_{N-2}^{t,-}$ be such that
$B'\preceq B$; here $t\in\ti0,t'\in\ti0$.

The following result was stated without proof in \cite{L22,4.4}.

(a) {\it If $t\ge0$, we have $t=t'$ or $|t'+1|<|t+1|$.  }
\nl
We have $\z(B')\le\z(B)$, see 2.10(c);
hence $\sha\z(B')^0\le\sha\z(B)^0$,
see 2.5(b). If $t'<0$, then $\sha\z(B)^0=\sha B^0=t$,
$\sha \z(B')^0=\sha B'{}^0=-t'-1$, so that $-t'-1\le t$; thus,
$|t'+1|<|t+1|$. If $t'\ge0$, then $\sha\z(B)^0=\sha B^0=t$,
$\sha\z(B')^0=\sha B'{}^0=t'$, so that $t'\le t$ and $t'+1\le t+1$.
This proves (a).

\subhead 3.5\endsubhead
The (fixed point free) involution of $E_N$ given by $X\m X^!=X+[1,N-1]$
leaves stable each of $E_N^+,E_N^-$. Let
$\bE_N,\bE_N^+,\bE_N^-$ be the
sets of orbits of this involution on
$E_N,E_N^+,E_N^-$ and let $\p:E_N@>>>\bE_N$,
$\p^+:E_N^+@>>>\bE_N^+$, $\p^-:E_N^-@>>>\bE_N^-$ be the orbit
maps.
Note that $X\m X^!$ interchanges $E_N^{t,+}$ with $E_N^{-t,+}$ (for any
$t\in\ti0,t\ne0$), ${}'E_N^{0,+}$ with ${}''E_N^{0,+}$ and
$E_N^{t,-}$ with $E_N^{-t-2,-}$ (for any $t\in\ti0$). Hence it
interchanges ${}'E_N^+$ with ${}''E_N^+$ and ${}'E_N^-$ with
${}''E_N^-$.

It follows that $\p^+$ restricts to a bijection ${}'E_N^+@>\si>>\bE_N^+$
and  $\p^-$ restricts to a bijection  ${}'E_N^-@>\si>>\bE_N^-$. 

For $Y\in\bE_N^+$ (resp. $Y\in\bE_N^-$) we denote by $\tY$ the unique
element of \lb
${}'E_N^+\cap(\p^+)\i(Y)$ (resp. ${}'E_N^-\cap(\p^-)\i(Y)$).

For $Y\in\bE_N$ we denote by $<Y>$ the image under
$\p$ of the subspace ${}'<\tY>$ of $E_N$. Note that $Y\in<Y>$.
If $Y\in\bE_N^+$ (resp. $Y\in\bE_N^+$) then
$\p^+$ (resp. $\p^-$) restricts to a bijection ${}'<\tY>@>>><Y>$
(resp. ${}'<\tY>\cap E_N^-@>\si>><Y>\cap\bE_N^-$); indeed, using
\cite{L23a, 5.4(a)}, we see that ${}'<\tY>$ (and even $<\tY>$)
does not contain
the kernel $F[1,N-1]$ of $\p:E_N@>>>\bE_N$.

\subhead 3.6\endsubhead
For $Y,Y'$ in $\bE_N^+$ we say that $Y'\bar\preceq Y$ if

(a) there exist $Y'=Y_0,Y_1,\do,Y_k=Y$ in $\bE_N^+$ such that
for any $i=0,1,\do,k-1$ we have either $\tY_i\preceq\tY_{i+1}$
or $|\g(\tY_i)|<|\g(\tY_{i+1})|$.

If for some $i$ we have $\tY_i\preceq\tY_{i+1}$, then
$|\g(\tY_i)|\le|\g(\tY_{i+1})|$. Indeed, if $B_i,B_{i+1}$ in
$\cx_{N-2}^+$ satisfy ${}'\e(B_i)=\tY_i$, ${}'\e(B_{i+1})=\tY_{i+1}$,
then $\sha B_i^0\le\sha B_{i+1}^0$; if $\sha B_i^0>0$, it follows that
$\g(\tY_i)<0,\g(\tY_{i+1})<0$ and $\g(\tY_i)-1\le-\g(\tY_{i+1})-1$
so that $|\g(\tY_i)|\le|\g(\tY_{i+1})|$; if $\sha B_i^0=0$, it follows that
$\g(\tY_i)=0$ and $|\g(\tY_i)|\le|\g(\tY_{i+1})|$ is obvious.

It follows that, if $Y'\bar\preceq Y$, then $|\g(Y')\le|\g(Y)|$;
moreover, if $|\g(Y')|=|\g(Y)|$ then (in the setup of (a)) we have
$\tY_i\preceq\tY_{i+1}$ for $i=0,1,\do,k-1$. We show:

(b) {\it  $\bar\preceq$ is a partial order on $\bE_N^+$.}
\nl
It is enough to show that if $Y\bar\preceq Y'$ and $Y'\bar\preceq Y$,
then $Y=Y'$. As we have seen above, we have $|\g(Y')|\le|\g(Y)|$
and similarly $|\g(Y)|\le |\g(Y')|$. It follows that
$|\g(Y')|=|\g(Y)|$.
Hence (in the setup of (a)) we have $\tY_i\preceq\tY_{i+1}$ for
$i=0,1,\do,k-1$, so that $\tY_0\preceq\tY_k$ that is $\tY'\preceq\tY$.
Similarly we have $\tY\preceq\tY'$. Since $\preceq$ is a partial order
it follows that $\tY=\tY'$ hence $Y'=Y$.

\subhead 3.7\endsubhead
For $Y,Y'$ in $\bE_N^-$ we say that $Y'\bar\preceq Y$ if

(a) there exist $Y'=Y_0,Y_1,\do,Y_k=Y$ in $\bE_N^-$ such that
for any $i=0,1,\do,k-1$ we have either $\tY_i\preceq\tY_{i+1}$ or
$|\g(\tY_i)+1|<|\g(\tY_{i+1})+1|$.
\nl
If for some $i$ we have $\tY_i\preceq\tY_{i+1}$ then
$|\g(\tY_i)+1|\le|\g(\tY_{i+1})+1|$. Indeed, if $B_i,B_{i+1}$ in
$\cx_{N-2}^-$ satisfy ${}'\e(B_i)=\tY_i$, ${}'\e(B_{i+1})=\tY_{i+1}$,
then $\sha B_i^0\le\sha B_{i+1}^0$; it follows that
$0\le\g(\tY_i)\le\g(\tY_{i+1})$, so that
$|\g(\tY_i)+1|\le|\g(\tY_{i+1})+1|$.

It follows that if $Y'\bar\preceq Y$, then
$|\g(Y')+1|\le |\g(Y)+1|$;
moreover, if $|\g(Y')+1|=|\g(Y)+1|$ then (in the setup of
(a)) we have
$\tY_i\preceq\tY_{i+1}$ for $i=0,1,\do,k-1$. We show:

(b) {\it $\bar\preceq$ is a partial order on $\bE^-_N$.}
\nl
It is enough to show that if $Y\bar\preceq Y'$ and $Y'\bar\preceq Y$,
then $Y=Y'$. As we have seen above we have $|\g(Y')+1|\le|\g(Y)+1|$
and similarly $|\g(Y)+1|\le|\g(Y')+1|$. It follows that
$|\g(Y')+1|=|\g(Y)+1|$. Hence (in the setup of (a)) we have
$\tY_i\preceq\tY_{i+1}$ for $i=0,1,\do,k-1$, so that $\tY_0\preceq\tY_k$,
that is $\tY'\preceq\tY$. Similarly we have $\tY\preceq\tY'$. Since
$\preceq$ is a partial order, it follows that $\tY=\tY'$ hence $Y=Y'$.

\subhead 3.8\endsubhead
Note that \cite{L22, 4.4-4.6} contains no proofs.
This and the next two subsections (3.9, 3.10) should be considered
as replacements (with proofs) of \cite{L22, 4.4-4.6}.

We show:

(a) {\it If $Y',Y$ in $\bE_N^+$ satisfy $Y'\in<Y>$, then
$Y'\bar\preceq Y$.}
\nl
Using our assumption we see that for some $X\in(\p_N^+)\i(Y')$
we have $X\preceq\tY$.
We have $(\p_N^+)\i(Y')=\{\tY',Z\}$ where $Z=(\tY')^!\in E_N^+$.
Assume first that $\g(\tY)<0$. If $X=\tY'$, then $\tY'\preceq\tY$ hence
$Y'\bar\preceq Y$. Thus we can assume that $X=Z\preceq\tY$. Since
$\g(\tY')\le0$, we have $\g(Z)\ge0$ so that $\g(Z)\ne\g(\tY)$; using
3.2(a) we see that $|\g(Z)|<|\g(\tY)|$ hence $|\g(\tY')|<|\g(\tY)|$.
We see that $Y'\bar\preceq Y$.

Next we assume that $\g(\tY)=0$ so that $\tY\in{}'E_N^{0,+}$.
If $X=\tY'$ then $\tY'\preceq\tY$, hence $Y'\bar\preceq Y$. Thus, we
can assume that $X=Z\in{}''E_N^{0,+}$. We have $Z\preceq\tY$. From 3.3(a)
we have $Z\in{}'E_N^{0,+}$, a contradiction. This proves (a).

Next we assume that $\g(\tY)=0$, so that $\tY\in{}'E_N^{0,+}$.
Then $Z\in{}''E_N^{0,+}$. Since $Z\preceq\tY$, from 3.3(a) we have
$Z\in{}'E_N^{0,+}$, a contradiction. This proves (b).

\subhead 3.9\endsubhead
We show:

(a) {\it If $Y',Y$ in $\bE_N^-$ satisfy $Y'\in<Y>$, then
$Y'\bar\preceq Y$.}
\nl
Using our assumption we see that for some $X\in(\p_N^-)\i(Y')$
we have $X\preceq\tY$.

We have $(\p_N^-)\i(Y')=\{\tY',Z\}$ where $Z=(\tY')^!\in E_N^-$,
$\g(Z)=-\g(\tY')-2$. We have $\g(\tY)\ge0$. If $X=\tY'$, then
$\tY'\preceq\tY$ hence $Y'\bar\preceq Y$. Thus we can assume that
$X=Z\preceq\tY$. Since $\g(\tY')\ge0$, we have $\g(Z)\le-2$ so that
$\g(Z)\ne\g(\tY)$; using 3.4(a) we see that
$|\g(Z)+1|<|\g(\tY)+1|$. We
have $\g(Z)+1=-\g(\tY')-1$, hence $|\g(\tY')+1|<|\g(\tY)+1|$. We see
that $Y'\bar\preceq Y$.

\subhead 3.10\endsubhead
Let $\ZZ[\bE_N^+]$ (resp. $\ZZ[\bE_N^-]$) be the abelian group consisting
of functions $\bE_N^+@>>>\ZZ$ (resp. $\bE_N^-@>>>\ZZ$).
This has a $\ZZ$-basis $\{\ps_Y;Y\in\bE_N^+\}$ (resp.
$\{\ps_Y;Y\in\bE_N^-\}$ where for $Y\in\bE_N^+$ (resp. $Y\in\bE_N^-$),
$\ps_Y\in\ZZ[\bE_N^+]$ (resp. $\ps_Y\in\ZZ[\bE_N^-]$) is the function
$Y\m1, Y'\m0$ if $Y'\in\bE_N^+-\{Y\}$ (resp. $Y'\in\bE_N^--\{Y\}$).
For any $Y\sub\bE_N^+$ (resp. $Y\sub\bE_N^-$) we set
let $\ti\ps_Y\in\ZZ[\bE_N^+]$ (resp. $\ti\ps_Y\in\ZZ[\bE_N^-]$)
be the characteristic function of $<Y>$ (resp.
$<Y>\cap\bE_N^-$).

From 3.8, 3.9, we see that
$\{\ti\ps_Y;Y\in\bE_N^+\}$ is related to $\{\ps_Y;Y\in\bE_N^+\}$ by an
upper triangular matrix with entries in $\ZZ$ and with $1$ on
diagonal
and that $\{\ti\ps_Y;Y\in\bE_N^-\}$ is related to
$\{\ps_Y;Y\in\bE_N^-\}$ by an upper triangular matrix with entries in
$\ZZ$ and with $1$ on diagonal. It follows that

$\{\ti\ps_Y;Y\in\bE_N^+\}$ is a $\ZZ$-basis of $\ZZ[\bE_N^+]$,

$\{\ti\ps_Y;Y\in\bE_N^-\}$ is a $\ZZ$-basis of $\ZZ[\bE_N^-]$.

\head 4. The equality ${}'\cx_{N-2}^+={}^\sha\cx_{N-2}^+$\endhead
\subhead 4.1\endsubhead
When $N\ge5$, for any $k\in[1,N-1]$ we define an (injective) map $\io_k:[1,N-2]@>>>[1,N]$ by
$\io_k(i)=i$ for $i\in[1,k-1]$, $\io_k(i)=i+2$ for $i\in[k,N-2]$. Note that
$\io_k(1)<\io_k(2)<\do<\io_k(N-2)$ and $\io_k(i)=i\mod2$ for all $i\in[1,N-2]$.
When $N\ge5$ and $k\in[1,N-1]$ we define a map $I_k:\cp_{N-2}@>>>\cp_N$ by
$$(i_1j_1,i_2j_2,\do,i_tj_t)\m(\io_k(i_1)\io_k(j_1),\io_k(i_2)\io_k(j_2),\do,\io_k(i_t)\io_k(j_t),
\{k,k+1\}).$$

(This is well defined since $\io_k:[1,N-2]@>>>[1,N]$ is injective with image not containing $k,k+1$.) 
\nl
If $N\ge5$, from the definitions, for $k\in[1,N-2]$, $t\in\ti0$ we have

(a) $I_k(\cx_{N-4}^{t,+})\sub\cx_{N-2}^{t,+}$,

(b) $I_k(\cx_{N-4}^{t,-})\sub\cx_{N-2}^{t,-}$

and for $k\in[1,N-3]$, we have

(c) $I_k({}'\cx_{N-4}^{0,+})\sub{}'\cx_{N-2}^{0,+}$

(d) $I_k({}''\cx_{N-4}^{0,+})\sub{}''\cx_{N-2}^{0,+}$. 

Let ${}'Pr_{N-2}^+$ be the subset of $\cp_N$ consisting of $\emp$
and of

$\{\{N-2,1\},\{N-3,2\},\do,\{N-\t,\t-1\}\}$
\nl
for various even $\t\in[2,(N-1)/2]$.

We define a subset ${}^\sha\cx_{N-2}^+$ of $\cp_N$ by induction on $N$ as
follows. We set ${}^\sha\cx_1^+=\emp$. Assume now that $N\ge5$.
Let $B\in\cp_N$. We say that $B\in{}^\sha\cx_{N-2}^+$ if either

$B\in{}'Pr_{N-2}^+$ or

$\sha B^0>0$ and there exists $B'\in{}^\sha\cx_{N-4}^+$ and $k\in[1,N-2]$
such that $B=I_k(B')$, or

$\sha B^0=0$ and there exists $B'\in{}^\sha\cx_{N-4}^+$ and $k\in[1,N-3]$
such that $B=I_k(B')$.

Note that ${}'Pr_{N-2}\sub{}'\cx_{N-2}^+$. Using (a),(c), we see
that when $N\ge5$ and $B'\in{}'\cx_{N-4}^+$ we have
$I_k(B')\in{}'\cx_{N-2}^+$ for $k\in[1,N-2]$ if $\sha B'{}^0>0$,
$I_k(B')\in{}'\cx_{N-2}^+$ for $k\in[1,N-3]$ if $\sha B'{}^0=0$.
From this we see by induction on $N$ that:

(e) {\it If $B\in{}^\sha\cx_{N-2}^+$ then $B\in{}'\cx_{N-2}^+$.}

\subhead 4.2\endsubhead
We show:

(a) {\it If $B\in{}'\cx_{N-2}^+$ then $B\in{}^\sha\cx_{N-2}^+$.}
\nl
We argue by induction on $N$. If $N=3$, (f) is obvious (both sides are
$\emp$). Now assume that $N\ge5$. By \cite{L22, 1.7(a)}, we have either

(b) $B\in Pr_{N-2}$ 
\nl
(notation of \cite{L22, 1.2}) or 

(c) there exists $B'\in\cx_{N-4}$ and $k\in[1,N-2]$ such that
$B=I_k(B')$.
\nl
Assume first that (b) holds. If $B\in Pr_{N-2}-{}'Pr_{N-2}^+$ then
from the definition we see that $B\n{}'\cx_{N-2}^+$, contradicting our
assumption. It follows that $B\in{}'Pr_{N-2}$ so that
$B\in{}^\sha\cx_{N-2}^+$.

Next we assume that (c) holds. If $B'\in\cx_{N-4}^-$, then from
4.1(b) we see that $B\in\cx_{N-2}^-$, contradicting our
assumption. Thus we have $B'\in\cx_{N-4}^+$.
If $B'\in\cx_{N-4}^{t,+}$ with $t>0$ then from 4.1(a) we see
that $B\in\cx_{N-2}^{t,+}$ with $t>0$, contradicting our assumption.
If $B'\in{}''\cx_{N-4}^{0,+}$ then from 4.1(d) we see that
$B\in{}''\cx_{N-2}^{0,+}$, contradicting our assumption.
If $B'\in{}'\cx_{N-4}^{0,+}$ and $k=N-2$ then $B\in{}''\cx_{N-4}^{0,+}$,
contradicting our assumption.
We see that we have either $B'\in\cx_{N-4}^{t,+}$ for some $t\in\un0,t<0$
or $B'\in{}'\cx_{N-4}^{0,+}$ and $k\in[1,N-3]$; the first alternative
occurs if $\sha B^0>0$ and the second alternative occurs if
$\sha B^0=0$.
Using the induction hypothesis we see that we have
$B'\in{}^\sha\cx_{N-4}^+$ and, in addition, if $\sha B^0=0$ we
have $k\in[1,N-3]$. It follows that $B\in{}^\sha\cx_{N-2}^+$,
proving (a).

Combining (a) with 4.1(e) we obtain the following result.

\proclaim{Theorem 4.3} We have ${}'\cx_{N-2}^+={}^\sha\cx_{N-2}^+$.
\endproclaim

\subhead 4.4\endsubhead
The set ${}^\sha\cx_{N-2}^+$ was introduced in \cite{L20, 2.1}. Several
properties of this set were stated without proof in the last paragraph
of \cite{L20, 2.2}. In view of 4.3, those properties are now seen to be
consequences of the corresponding properties of ${}'\cx_{N-2}^+$, which
are proved in this paper.

\subhead 4.5  \endsubhead
Let ${}'Pr_{N-2}^-$ be the subset of $\cp_N$ consisting of
$\{N-1,N\}$ and of

$\{\{N-1,N\},\{N-2,1\},\{N-3,2\},\do,\{N-1-\t,\t\}\}$
\nl
with $\t$ even, $\t\in[2,(N-3)/2]$.

We define a subset ${}^\sha\cx_{N-2}^-$ of $\cp_N$ by induction on $N$ as
follows. We set ${}^\sha\cx_1^-=\{N-1,N\}$. Assume now that $N\ge5$.
Let $B\in\cp_N$. We say that $B\in{}^\sha\cx_{N-2}^-$ if either

$B\in{}'Pr_{N-2}^-$ or

there exists $B'\in{}^\sha\cx_{N-4}^-$ and $k\in[1,N-2]$
such that $B=I_k(B')$.
\nl
Using 4.1(b), we seethat when $N\ge5$ and $B'\in{}'\cx_{N-4}^-$ we have
$I_k(B')\in{}'\cx_{N-2}^+$ for $k\in[1,N-2]$.
From this we see by induction on $N$ that:

(a) {\it If $B\in{}^\sha\cx_{N-2}^-$ then $B\in{}'\cx_{N-2}^-$.}

\subhead 4.6  \endsubhead
We show:

(a) {\it If $B\in{}'\cx_{N-2}^-$ then $B\in{}^\sha\cx_{N-2}^-$.}
\nl
We argue by induction on $N$. If $N=3$, (a) is obvious (both sides are
$\{N-1,N\}$). Now assume that $N\ge5$. By \cite{L22, 1.7(a)}, we have
either

(b) $B\in Pr_{N-2}$ 
\nl
(notation of \cite{L22, 1.2}) or

(c) there exists $B'\in\cx_{N-4}$ and $k\in[1,N-2]$
such that $B=I_k(B')$.

Assume first that (b) holds. If $B\in Pr_{N-2}-{}'Pr_{N-2}^-$ then
from the definition we see that $B\n{}'\cx_{N-2}^-$, contradicting our
assumption. It follows that $B\in{}'Pr_{N-2}^-$ so that
$B\in{}^\sha\cx_{N-2}^-$.

Next we assume that (c) holds.
If $B'\in\cx_{N-4}^+$ then from 4.1(a) we see that
$B\in\cx_{N-2}^+$, contradicting our assumption.
Thus we have $B'\in\cx_{N-4}^-$.
If $B'\in\cx_{N-4}^{t,-}$ with $t<0$ then from 4.1(b) we see
that $B\in\cx_{N-2}^{t,+}$ with $t<0$, contradicting our assumption.
We see that we have $B'\in\cx_{N-4}^{t,-}$ for some $t\in\un0,t\ge0$.
Using the induction hypothesis we see that we have
$B'\in{}^\sha\cx_{N-4}^-$. It follows that $B\in{}^\sha\cx_{N-2}^-$, proving
(a).

Combining (a) with 4.5(a) we obtain the following result.

\proclaim{Theorem 4.7} We have ${}'\cx_{N-2}^-={}^\sha\cx_{N-2}^-$.
\endproclaim

\head 5. Fourier transform\endhead
\subhead 5.1\endsubhead
For $X,X'$ in $E_N$ we set $(X,X')=\sha(X\cap X')\mod2$. This is a
nondegenerate symplectic form on $E_N$ with values in $F$.
Let $\QQ[E_N]$ be the $\QQ$-vector space of functions $E_N@>>>\QQ$.
Recall that the Fourier transform $\Ph:\QQ[E_N]@>>>\QQ[E_N]$
is the involution $f\m\Ph(f)$ where

$\Ph(f)(X)=2^{-(N-1)/2}\sum_{X'\in E_N}(-1)^{(X,X')}f(X')$
\nl
for $X\in E_N$. 
Recall that $\QQ[E_N]$ has a $\QQ$-basis
$\{\ps_{<X>};X\in E_N\}$; applying $\Ph$ to this basis we obtain
another basis $\{\Ph(\ps_{<X>});X\in E_N\}$ of $\QQ[E_N]$.
For $X\in E_N$ we have 
$\Ph(\ps_{<X>})=\sum_{X'\in E_N}A_{X',X} \ps_{<X'>}$ where
$A_{X',X}\in\QQ$. According to \cite{L20a}, the matrix
$(A_{X',X})$ is upper triangular for a suitable order of the set
of indices. More precisely, if for $k\ge0$ we denote by
$\QQ[E_N]^{\ge k}$ the subspace of $\QQ[E_N]$ spanned by all
$\ps_{<X>}$ with $X\in E_N,\dim<X>\ge k$, we have that
$\QQ[E_N]^{\ge k}$ is stable under $\Ph$ and that the linear map
induced by $\Ph$ on $\QQ[E_N]^{\ge k}/\QQ[E_N]^{\ge k+1}$ is
$\pm$ the identity.

Similarly, applying $\Ph$ to the basis $\{\ps_{{}'<X>};X\in E_N\}$; we
obtain another basis $\{\Ph(\ps_{{}'<X>});X\in E_N\}$ of $\QQ[E_N]$.
For $X\in E_N$ we have 
$$\Ph(\ps_{{}'<X>})=\sum_{X'\in E_N}{}'A_{X',X}\ps_{{}'<X'>}$$ where
${}'A_{X',X}\in\QQ$. We show that the matrix
$({}'A_{X',X})$ is in general not upper triangular for any order
of the set of indices. Assume that $N=5$.
We have $\QQ[E_N]=\QQ[E_N]_2\op\QQ[E_N]_1\op\QQ[E_N]'_1\op\QQ[E_N]_0$
where $\QQ[E_N]_2$ is spanned by
all $\ps_{{}'<X>}$ where $X\in E_N$ satisfies $\dim<X>=2$;
$\QQ[E_N]_1$ is spanned by all $\ps_{{}'<X>}$ where $X\in E_N$ satisfies
$\dim<X>=1$ and $X$ is not $45,51,42,13$;
$\QQ[E_N]'_1$ is spanned by all $\ps_{{}'<X>}$ where $X\in E_N$
is one of $45,51,42,13$;
$\QQ[E_N]'_0$ is spanned by  $\ps_{{}'<\emp>}$.
Then the subspaces

$\QQ[E_N]_2$, $\QQ[E_N]_2\op\QQ[E_N]_1$,
$\QQ[E_N]_2\op\QQ[E_N]_1\op\QQ[E_N]'_1$
\nl
are stable under $\Ph$. By passage to quotients, $\Ph$ induces on
$\QQ[E_N]_2$, $\QQ[E_N]_1$, $\QQ[E_N]'_1$, $\QQ[E_N]_0$, linear maps
denoted by $\a_2,\a_1,\a'_1,\a_0$ respectively. We have
$\a_2=1,\a_1=-1,\a_0=-1$. Moreover, $2\a'_1$ is given by the matrix
$$\left(\matrix-1&-1&0&1\\-1&-1&1&0\\-2&2&1&1\\2&-2&1&1\endmatrix
\right)$$

which is not upper triangular for any order of the index set.

\subhead 5.2\endsubhead
$(,)$ restricts to a symplectic form on $E_N^+$ with one dimensional
radical spanned by $[1,N-1]$ and this induces a nondegenerate
symplectic form $(,)^+$ on $\bE_N^+$.
Let $\QQ[\bar E_N^+]$ be the $\QQ$-vector space of functions
$\bE_N^+@>>>\QQ$. Then the Fourier transform
$\bar\Ph:\QQ[\bE_N^+]@>>>\QQ[\bE_N^+]$ is defined in terms of $(,)^+$.
It is likely that the matrix of this linear map with respect to
the basis

(a) $\{\ti\ps_Y;Y\in\bE_N^+\}$
\nl
of $\QQ[\bE_N^+]$ is upper
triangular for a suitable partial order of the set of indices.
This is indeed so if $N\le7$.

\enddocument